\ifx\documentclass\undefined
\ifx\selectfont\undefined
\documentstyle[11pt]{article}
\else
\documentstyle[oldlfont,11pt]{article}
\fi
\else
\documentclass[11pt]{article}
\usepackage{latexsym}
\fi
\font\tenmsb=msbm10  at 11pt
\font\sevenmsb=msbm7 at 9pt
\font\fivemsb=msbm5 at 7pt
\newfam\msbfam
\textfont\msbfam=\tenmsb
\scriptfont\msbfam=\sevenmsb
\scriptscriptfont\msbfam=\fivemsb
\def\Bbb#1{\fam\msbfam\relax#1}
        \def\nat{{\Bbb N}}
        \def\real{{\Bbb R}}
        
        \def\F{{\cal F}}

        \def\T{{\cal T}}
\def\supp{\mathop{\rm supp}\nolimits}
\def\Tr{\mathop{\rm Tr}\nolimits}
\def\ep{\varepsilon}

\def\To{\Rightarrow}
\def\frac#1#2{{#1\over#2}}
\def\EQ#1{(\ref{eq:#1})}
\newtheorem{fact}{Fact}[section]
\newtheorem{thm}[fact]{Theorem}
\newtheorem{lemma}[fact]{Lemma}
\newtheorem{prop}[fact]{Proposition}
\newtheorem{cor}[fact]{Corollary}

\newtheorem{remark}[fact]{Remark}

\newtheorem{problem}[fact]{Problem}
\newenvironment{proof}{\medskip \par \noindent {\it Proof.}\ }{\hfill
	$\Box$ \medskip \par}
\def\pf#1{\medskip \par \noindent {\it #1.}\ }
\def\endpf{\hfill $\Box$ \medskip \par}
\begin{document}
\title{On certain equivalent norms on Tsirelson's space}
\author{Edward  Odell \thanks{Research supported by
NSF and TARP.}
\and Nicole Tomczak-Jaegermann\thanks{This author held the Canada
  Council Killam Research Fellowship in 1997/99.}}
\date{}
\maketitle
\begin{abstract}
Tsirelson's space $T$ is known to be distortable but it is open as to 
whether or not $T$ is arbitrarily distortable. 
For $n\in \nat$ the norm $\|\cdot\|_n$ of the Tsirelson space 
$T(S_n,2^{-n})$ is equivalent to the standard norm on $T$. 
We prove there exists $K<\infty$ so that for all $n$, $\|\cdot\|_n$ 
does not $K$ distort any subspace $Y$ of $T$.
\end{abstract}

\section{Introduction}\label{sec0}
\setcounter{equation}{0}

It remains open an important question as to whether or not there
exists a distortable Banach space which is not arbitrarily
distortable.  The primary candidate for such a space is Tsirelson's
space $T$.  While it is not difficult to directly define, for every $1
< \lambda <2$, an equivalent norm on $T$ which is a
$\lambda$-distortion, $T$ does not belong to any general class of
Banach spaces known to be arbitrarily distortable.  In fact (see
below) if there does exists a distortable not arbitrarily distortable
Banach space $X$ then $X$ must contain a subspace which is very
Tsirelson-like in appearance.  It is thus of interest to examine in
particular all known equivalent norms on $T$ to see if they can
arbitrarily distort $T$ (or a subspace of $T$). We do so in this paper
for a previously unstudied fascinating class of renormings.

The renormings we consider here are ``natural'' in that pertain to the
deep combinatorial nature of the norm of $T$.  Namely, for each $n$ by
$\|\cdot\|_n$ we denote the norm of the Tsirelson space
$T(S_n,2^{-n})$, which can be easily seen to be equivalent to the
original norm on $T$.  Our main result (Theorem~\ref{thm:2.1}) is that
this family of equivalent norms does not arbitrarily distort $T$ or
even any subspace of $T$.  The proof actually introduces a larger
family of equivalent norms $(\|\cdot\|_j^n)_{j,n}$ and
$(|\cdot|_j^n)_{j,n}$ which are shown to not arbitrarily distort any
subspace of $T$.  Quantitative estimates for the stabilizations of
these norms are given in Theorem~\ref{thm:3.1}.  It is shown that (up
to absolute constants) one has that for all $n$ and subspaces $X
\subseteq T$ there is a subspace $Y \subseteq X$ such that $\|y\|_n
\sim \frac{1}{n}$ if $y \in Y$ with $\|y|| =1$.

Some stabilization results for more general norms on $T$ of various
classes are also given in Section~\ref{sec3}.  In Section~\ref{sec4}
we raise some problems.

Section~\ref{sec1} contains the relevant terminology and background
material.  Otherwise our notation is standard as may be found in
\cite{LT}.

More detailed information about Tsirelson's space and Tsirelson type
spaces can be found in \cite{CS}, \cite{OTW}, \cite{AD}, \cite{AO} and
the references therein.  

\section{Preliminaries}\label{sec1}
\setcounter{equation}{0}

$X,Y,Z,\ldots$ will denote separable infinite dimensional real Banach
spaces.
If $(x_i)$ is a basic sequence, $(y_i) \prec (x_i)$ shall mean that $(y_i)$ 
is a block basis of $(x_i)$. 
$X= [(x_i)]$ is the closed linear span of $(x_i)$. 
If $X$ has a basis $(x_i)$, $Y\prec X$ denotes $Y= [(y_i)]$ where 
$(y_i)\prec (x_i)$. 
The terminology is imprecise in that ``$\prec$'' refers to a fixed basis 
for $X$ but no confusion shall arise. 
 $S_X = \{x\in X: \|x\|= 1\}$.

A space $(X,\|\cdot\|)$ is  {\it arbitrarily distortable}, 
if for all $\lambda >1$ there exists an equivalent norm $|\cdot|$ on $X$ 
satisfying for all $Y\subseteq X$ 
\begin{equation}\label{eq:1.1}
\sup \left\{ {|y|\over |z|} : y,z\in S_Y\right\} >\lambda\ .
\end{equation} 
The norm $|\cdot|$ satisfying \EQ{1.1} is said to $\lambda$-distort $X$. 
$X$ is {\it $\lambda$-distortable\/} if some norm $\lambda$-distorts $X$.
$X$ is {\it distortable\/} if it is $\lambda$-distortable for some 
$\lambda >1$. 
If $X$ has a basis then ``for all $Y\subseteq X$'' in the statement
above can be replaced by  ``for all $Y\prec X$''.

Tsirelson's space $T$ (defined below) is known to be $2-\ep$
distortable for all $\ep>0$ (see e.g., \cite{OTW}).  If a space
$X$ exists which is distortable but not arbitrarily distortable then
$X$ can be assumed to have an unconditional basis \cite{T-J}, to
be asymptotic $c_0$ or $\ell_p$ for some $1\le p<\infty$ \cite{MT}
and to contain $\ell_1^n$'s uniformly \cite{M}. These
characteristics in conjunction with others developed in \cite{OTW}
yield that $T$ is the prime candidate for such a space.

For $n\in \nat$, the Schreier classe $S_n$ is a pointwise compact
hereditary collection of finite subsets of $\nat$ \cite{AA}.
We write for $E,F\subseteq \nat$, $E<F$ (resp., $E \le F$) if $\max
E<\min F$ (resp., $\max E\le \min F$) or if either one is empty.
$$S_0= \Big\{ \{n\} :n\in \nat \Big\} \cup \{\emptyset\}\ .$$ 
We inductively define 
$$S_{k+1} = \biggl\{ \bigcup_{p=1}^\ell E_p : \{\ell\} \le E_1 < \cdots <
E_\ell \mbox{ and } E_p \in S_k\mbox{ for } 1\le p\le \ell\biggr\}\ .$$
$(E_i)_{i=1}^\ell$ is $k$-admissible if $E_1 <\cdots < E_\ell$ and 
$(\min E_i)_{i=1}^\ell \in S_k$. 
It is easy to see that 
\begin{eqnarray*}
S_k [S_n] & \equiv & \biggl\{ \bigcup_{i=1}^\ell E_i : (E_i)_1^\ell 
\mbox{ is $k$-admissible and } E_i \in S_n \mbox{ for } 1 \le i\le
\ell  \biggr\}\\
&=& S_{n+k}\ .
\end{eqnarray*}

If $(y_i)$ is a basis then $(x_i)_1^\ell \prec (y_i)$ is 
{\it $k$-admissible\/} (w.r.t.\ $(y_i)$) if $(\supp x_i)_{i=1}^\ell$ is 
$k$-admissible. 
Here if $x= \sum_{i\in A} a_i y_i$ and $a_i\ne0$ for $i\in A$, then 
$\supp x=A$. 

$c_{00}$ denotes the linear space of finitely supported real sequences and 
$(e_i)$ is the unit vector basis for $c_{00}$. 
If $x = \sum_i x(i) e_i \in c_{00}$ and $E\subseteq \nat$ then $Ex \in
c_{00}$ is defined by $Ex= \sum_{i\in E} x(i) e_i$.
Let $\F$ be a pointwise compact {\it hereditary\/} 
(that is, $G\subseteq F\in \F \To G\in \F$) 
family of finite subsets of $\nat$ containing $S_0$ and let $0<\lambda<1$. 
The Tsirelson space $T(\F,\lambda)$ is the completion of $c_{00}$ under 
the implicit norm 
\begin{equation}\label{eq:1.2} 
\|x\| = \|x\|_\infty \vee \sup \biggl\{ \lambda \sum_{i=1}^\ell \|E_ix\| 
: E_1 < \cdots < E_\ell \mbox{ and } (\min E_i)_1^\ell \in \F\biggr\}
\end{equation}
$(e_i)$ is then a normalized unconditional basis for $T(\F,\lambda)$. 
Furthermore if $\F \supseteq S_1$ then $T(\F,\lambda)$ does not
contain an isomorph of $\ell_1$ but is {\it asymptotically\/} $\ell_1$
(that is, if $(x_i)_1^\ell$ is 1-admissible then $\|\sum_1^\ell x_i\|
\ge \lambda \sum_1^\ell \|x_i\|$).
The existence of such a norm \EQ{1.2}  can be found in \cite{AD}.

The classical Tsirelson's space is $T\equiv T(S_1,2^{-1})$ and we
write $\|\cdot\|$ ($ = \|\cdot\|_1$) for the norm of $T$.
We also consider the space $T(S_n,2^{-n})$, for a fixed $n \in \nat$,
and we denote its norm by $\|\cdot\|_n$.
These norms are all equivalent on $c_{00}$ and thus the spaces
coincide.
Indeed,
\begin{equation}\label{eq:1.3} 
\|x\|_n \le \|x\| \le 2^{n-1} \|x\|_n\mbox{ for } x\in T\ .
\end{equation} 

We explain \EQ{1.3} and set some terminology for later use. 
$\|x\|$ is calculated as follows. 
If $\|x\| \ne \|x\|_\infty$ then $\|x\| = \frac12 \sum_1^\ell\|E^1_ix\|$ for 
some 1-admissible collection $(E^1_i)_1^\ell$. 
For $i\le \ell$ either $\|E^1_ix\| = \|E^1_i x\|_\infty$ or $\|E^1_i x\|$ 
is calculated by means of a similar decomposition. 
Ultimately one obtains for some finite $A\subseteq \nat$,  
$$\|x\| = \sum_{i\in A} 2^{-n(i)} |x(i)|,$$
where $n(i)$ is the number of decompositions necessary before obtaining a 
set $E_j^{n(i)}$ for which $\|E_j^{n(i)}x\| = \|E_j^{n(i)} x\|_\infty 
= |x(i)|$. 

Thus the norm in $T$ can be described as follows in terms of trees of sets. 
By an {\it admissible tree $\T$ of sets\/} we shall mean $\T = (E_i^n)$ 
for $1\le i\le i(n)$, $0\le n\le k$  is a tree of finite subsets of 
$\nat$ partially ordered by reverse inclusion with the following properties. 
$E_i^n$ is said to have level $n$. 
$i(0)=1$, $E_i^n < E_j^n $ if $i< j$, all successors of any $E_i^n$
form a 1-admissible partition of $E_i^n$ and every set $E_i^{n+1}$ is
a successor of some $E_j^n$.
Thus all sets of level $n$ form an $n$-admissible collection.  $E_i^n$
is a {\it terminal\/} set of $\T$ if it has no successors.

Thus one has for $x\in T$ 
\begin{eqnarray}\label{eq:1.4}
\|x\| &=& \sup \biggl\{ \sum_{i\in A} 2^{-n(i)} \|E_ix\|_\infty : 
(E_i)_{i\in A} \mbox{ are terminal sets}\\
&&\mbox{of an admissible tree with level } E_i = n(i)\biggr\}\ .\nonumber 
\end{eqnarray}

Also \EQ{1.4} holds if $\|E_ix\|_\infty$ is replaced by $\|E_ix\|$. 

The norm $\|\cdot\|_n$ is calculated in a similar fashion except that
terminal sets are allowed only to have levels $kn$ for some
$k=0,1,2,\ldots$
\begin{eqnarray}\label{eq:1.5} 
\|x\|_n &=& \sup \biggl\{ \sum_{i\in A} 2^{-nk(i)} \|E_ix\|_\infty : 
(E_i)_{i\in a} \mbox{ are terminal sets of an admissible tree}\\
&&\mbox{where  $E_i$ has level $nk(i)$ for some } k(i) = 0,1,2,\ldots\biggr\}
\ .\nonumber
\end{eqnarray}

{From} these formulas we see that $\|x\|_n\le \|x\|$. 
Furthermore if $\T$ is an admissible tree, terminal sets not having
levels $0,n,2n,\cdots$ can be continued to the next such level, an
increase of at most $n-1$ levels, yielding $\|x\| \le 2^{n-1}
\|x\|_n$.

More exotic {\it mixed Tsirelson spaces\/} were introduced in
\cite{AD}.
We shall not discuss a general definition, but we shall give a formula
for the norm in a special case of interest here.
For $j\ge0$ and $n\in \nat$ we let $\|\cdot\|_j^n$ be the norm of the mixed 
Tsirelson space $T\left((S_{j+kn},2^{-(j+kn)})_{k=0}^\infty\right)$. 
One obtains a similar formula to that in \EQ{1.4} for the norm except that 
terminal sets may only have levels $j,j+n,j+2n,\ldots$ 
\begin{eqnarray}\label{eq:1.6} 
\|x\|_j^n & = & \|x\|_\infty \vee \sup \biggl\{ \sum_{i\in A} 2^{-(j+nk(i))} 
\|E_i x\|_\infty : (E_i)_{i\in A} \mbox{ are terminal sets}\\
&&\mbox{of an admissible tree having level } E_i = j+nk(i) \mbox{ for some }
k(i)= 0,1,2,\ldots\biggr\}\ .\nonumber
\end{eqnarray}
Thus $\|\cdot\|_0^n = \|\cdot\|_n$. 
Furthermore $\|\cdot\|_j^n$ is an equivalent norm on $T$. 

We prove in Section 2 that the family of norms $(\|\cdot\|_j^n)_{n,j}$ cannot 
arbitrarily distort any subspace of $T$. 
We do this by introducing a slight variation  of $\|\cdot\|_j^n$
(which omits the first term in \EQ{1.6}):
\begin{eqnarray}\label{eq:1.7} 
|x|_j^n & = & \sup \biggl\{ \sum_{i\in A} 2^{-(j+nk(i))} \|E_ix\|_\infty : 
(E_i)_{i\in A} \mbox{ are  terminal sets}\\
&&\mbox{of an admissible tree having level } E_i = j+nk(i),\ k(i)\ge 0\biggr\}
\ .\nonumber
\end{eqnarray}
Thus $|\cdot|_0^n = \|\cdot\|_n$, $|\cdot|_j^n$ is an equivalent norm
on $T$ an $|\cdot|_j^n \le \|\cdot\|_j^n $.
Our next proposition  yields some simple facts about $|\cdot|_j^n$. 
Statements a) and b) are the reason we work with $|\cdot|_j^n$ rather than 
directly with $\|\cdot\|_j^n$. 
Moreover d) yields 
that $\|\cdot\|_j^n$ and $|\cdot|_j^n$ are nearly 
the same on some subspace of any given $Y\prec T$. 
First recall the Schreier space $X_m$ (\cite{AA}, also
\cite{CS}, for $m =1$).
$X_m$  is the completion of $c_{00}$ under 
$$|x|_m = \sup \biggl\{ \Big| \sum_{i\in E} x(i)\Big| : E\in S_m\biggr\}\ .$$
$X_m$ is isometric to a subspace of $C(\omega^{\omega^m})$ and hence is 
$c_0$-saturated: if $Y\subseteq X$ then $Y$ contains an isomorph of $c_0$. 
For $Z\subseteq T$, $S_Z$ is the unit sphere w.r.t.\ the Tsirelson norm 
$\|\cdot\|$. 

\begin{prop}\label{prop:1.1}
{\rm a)} Let $j\ge0$ and $n\in \nat$. For $x\in T$ 
$$|x|_j^n = \frac1{2^j} \sup \biggl\{ \sum_{\ell=1}^r \|E_\ell x\|_n 
: (E_\ell)_1^r\mbox{ is $j$-admissible}\biggr\}$$

{\rm b)} Let $j\ge0$ and $k,n\in \nat$. For $x\in T$
$$|x|_{j+k}^n = \frac1{2^k} \sup \biggl\{ \sum_{\ell=1}^r |E_\ell x|_j^n : 
(E_\ell)_1^r \mbox{ is $k$-admissible}\biggr\}$$ 

{\rm c)} Let $\ep >0$, $n,k\in \nat$ and $0\le j< n$. 
Let $Y\prec T$. 
Then there exists $Z\prec Y$ so that if $z\in S_Z$
\begin{equation}\label{eq:1.8} 
\Big|\, |z|_j^n - |z|_{j+np}^n\,\Big| < \ep\mbox{ if } 1\le p\le k\ .
\end{equation}

{\rm d)} For $n,j\in \nat$, $\ep>0$ and $Y\prec T$ there exists $Z\prec Y$ 
so that for all $z\in S_Z$, 
$$\Big|\, |z|_j^n - \|z\,\|_j^n\Big| <\ep\mbox{ and } 
\Big|\, |z|_n^n - \|z\|_n\,\Big| <\ep\ .$$ 
\end{prop}

\begin{proof}
a) and b) follow easily from \EQ{1.5}--\EQ{1.7} and the fact that 
$S_{k+j} = S_k[S_j]$. 
c) is proved by choosing $Z$ so that the first few levels of the admissible 
tree used to compute $|z|_{j+nk}^n$ will contribute only a negligible amount. 
Precisely we first note that 
$$|z|_j^n \ge |z|_{j+np}^n \ge |z|_{j+nk}^n\mbox{ for } 1\le p\le k\ .$$
Thus we need only achieve \EQ{1.8} for $p=k$. 
Let $|\cdot|_{j+nk}$ be the norm of the Schreier space $X_{j+nk}$. 
For $z\in T$ let 
$$|z|_j^n = \sum_{\ell\in A} 2^{-(j+nk(\ell))} |z(\ell)|$$ 
be obtained from  \EQ{1.7}. 
Thus if 
$$E= \{\ell\in A: k(\ell) <k\}$$
then $E\in S_{j+nk}$ and so 
$$|z|_j^n \le |Ez|_{j+nk} + |z|_{j+k}^n \le |z|_{j+nk} + |z|_{j+k}^n\ .$$ 
Also $|z|_{j+nk} \le 2^{j+nk} \|z\|$ for $z\in  T$. 
Since $X_{j+nk}$ is $c_0$-saturated and $T$ does not contain $c_0$ it follows 
that given $Y\prec T$ there exists $Z\prec Y$ so that if $z\in S_Z$ then 
$|z|_{j+nk}<\ep$. 
This proves c). 
d) is proved similarly to c). 
The norms in question differ only in that the terminal sets of an admissible 
tree can differ only in a finite number of levels. 
\end{proof}

We shall need a generalized notion of $n$ admissible. 
For $k,n\in\nat$, 
$(E_r)_1^s$ is {\it $n$ admissible\/} $(k)$ if $(kE_r)_1^s$ is $n$ admissible 
where $kE \equiv \{ke:e\in E\}$. 
Similarly we say $(y_i)_1^\ell \prec (e_i)$ is $n$ admissible $(k)$ if 
$(\supp y_i)_1^\ell$ is $n$ admissible $(k)$. 
Also we say a tree $\T$ is admissible $(k)$ if $(kE)_{E\in \T}$ is an 
admissible tree.

\begin{prop}\label{prop:1.2}
There exists $K_1<\infty$ so that if $n,k\in \nat$, $1>\ep >0$ and 
$(y_i) \prec (e_i)$ is normalized (in $T$), then there exists a finite set 
$A\subseteq \nat$ and $(\alpha_\ell)_{\ell\in A} \subset (0,1]$ so that 
$(y_\ell)_{\ell\in A}$ is $n$ admissible and setting 
$z= \sum_{\ell\in A} \alpha_\ell y_\ell$ we have the following: 

{\rm i)}\ \ \ $\sum_{\ell\in A} \alpha_\ell = 2^n$. 

{\rm ii)}\ \ If $B\subseteq A$ and $(y_\ell)_{\ell\in B}$ 
is $n-1$ admissible $(k)$ then $\sum_{i\in B}\alpha_i <\ep$.

{\rm iii)} $1\le \|z\|\le K_1$
\end{prop}

We call such a 
$z$ an {\it $(n,\ep)$ average\/} $(k)$ of $(y_\ell)$. 
This was proved in \cite{OTW} for $k=1$. 
The proof uses the following fact (see e.e., \cite{CS}, Prop.~II.4).

\begin{prop}\label{prop:1.3} 
 There exists $K_2<\infty$ so that if $(y_i)$ is a normalized block 
basis of $(e_i)$ in $T$ then for all $(a_i)$ if 
$m_i = \min \supp y_i$, 
$$\| \sum a_i e_{m_i}\| \le \|\sum a_i y_i\| 
\le K_2 \|\sum a_i e_{m_i}\|\ .$$
\end{prop}

\pf{Proof of Proposition~\ref{prop:1.2}} 
By passing to a subsequence of $(y_i)$ we may assume that $m_{i+1}>km_i$ 
where $m_i = \min \supp y_i$. 
By \cite{OTW} we can find $z= \sum_{\ell\in A} \alpha_\ell y_\ell$, 
$(\alpha_\ell)_{\ell\in A} \subseteq \real^+$, 
$\sum_{\ell\in A} \alpha_\ell = 2^n$ and  
$\sum_{\ell\in B} \alpha_\ell < \ep/2$ if 
$(m_\ell)_{\ell\in B} \in S_{n-1}$. 
Furthermore $1\le \|z\| \le K_1$. 
It remains to check that ii) holds. 
Suppose that $B\subseteq A$ so that $(km_i)_{i\in B} \in S_{n-1}$. 
Since $m_{i+1} >km_i$ this yields that $(m_{i+1})_{i\in B} \in S_{n-1}$ 
and hence $(m_i)_{i\in B\setminus \min B}\in S_{n-1}$. 
Thus $\sum_{\ell\in B\setminus\min B}\alpha_\ell <\ep/2$. 
Also $\alpha_{\min B} <\ep/2$ and so ii) holds.
\endpf

\section{Stabilizing the norms $(\|\cdot\|_n)$}\label{sec2}
\setcounter{equation}{0}

Our goal is to prove that the norms $(\|\cdot\|_j^n)$ and  hence in
particular the norms $(\|\cdot\|_n)$ do not arbitrarily distort any
subspace of $T$.
In light of Proposition~\ref{prop:1.1} it suffices to prove 

\begin{thm}\label{thm:2.1}
There exists $K>1$ so that for all $Y\prec T$ and $n\in \nat$ there exist
$Z\prec Y$ and $d>0$ satisfying: for all $0\le j< n$ and 
$z\in S_Z$ 
$$d\le |z|_j^n \le Kd$$
\end{thm} 

Before beginning the proof we recall that there exists $K_3< \infty$
so that $\|\sum b_i e_{3i}\| \le K_3 \|\sum b_i e_i\|$ \cite{CS},
Prop.~I.12.

\begin{lemma}\label{lem:2.2}
Let $(w_i)$ be a normalized block basis of $(e_i)$ in $T$. 
Suppose that for some $c>0$ and $L\ge 1$, we have for all $i$ 
$$L^{-1} c\le |w_i|_j^n \le Lc\ \mbox{ for }\ 0\le j\le n\ .$$ 
Let $w = \sum a_i w_i$, $\|w\|=1$. 
Then for $0\le j< n$, $c(LK_2)^{-1} \le |w|_j^n \le 2LK_3 c$.
\end{lemma}

\begin{proof} 
{From} Proposition~\ref{prop:1.3} there exists an admissible tree $\T$
whose terminal sets are all equal to $\supp w_i$ for some $i$ yielding
$$\|w\| \ge \sum_{i\in A} |a_i| 2^{-n(i)} \|w_i\|
= \sum_{i\in A} |a_i| 2^{-n(i)} \ge K_2^{-1}\ .$$
Let $1\le j\le n$ be fixed. 
We shall produce a lower estimate for $|w|_j^n$ by extending $\T$ as 
follows. 
Fix $i\in A$ and consider the term $|a_i| 2^{-n(i)} \|w_i\|$. 
Suppose this term resulted from $E=\supp w_i$ where $E$ was terminal in $\T$ 
of level $n(i)$.
First suppose that $n(i) \ge j$ so that $n(i) = j+kn+p$ for some $0\le
p< n$ and $k\ge0$; then let $q=n-p$. 
If $n(i) <j$ let $q= j-n(i)$. 
If $q\ge1$ extend $\T$ $q$-levels below $E$ via the $q$-admissible 
family of sets which yield by Proposition~\ref{prop:1.1}
$$|w_i|_q^n = \frac1{2^q} \sum_{s=1}^r \|E_s^i w_i \|_n \ge cL^{-1}\ .$$
The new tree has terminal sets only at levels $(j+kn)_{k=0}^\infty$. 
When used in \EQ{1.7} it yields 
$$|w|_j^n \ge \sum |a_i| 2^{-n(i)} cL^{-1} \ge c(LK_2)^{-1}\ .$$

For the upper estimate let $\T$ be the admissible tree having terminal sets 
(which we may assume to be singletons) 
of levels $j,j+n,j+2n,\ldots$ which produces $|w|_j^n$ in \EQ{1.7}. 
We say $w_i$ is {\it badly split\/} by some level of $\T$ if there exist 
$E\ne F$ in $\T$ having the same level with $Ew_i\ne0$, $Ew_s\ne0$ for some 
$s\ne i$ and $Fw_i \ne0$. 
If no $w_i$ is badly split by some level of $\T$ then if for some $i$,
$\supp w_i$ contains a terminal set in $\T$ there exists a
1-admissible family $(E_s^i)_1^{\ell (i)}$ in $\T$ of minimal level
having the property that $\bigcup_1^{\ell(i)} E_s^i\subseteq \supp
w_i$ and $F\cap \supp w_i=\emptyset$ for all other $F\in \T$ of the
same level as the $E_s^i$'s.
Thus for some set $A$, 
\begin{equation}\label{eq:2.1} 
|w|_j^n = \sum_{i\in A} 2^{-n(i)} |a_i| \sum_{s=1}^{\ell (i)} 
|E_s^i w_i|_{j(i)}^n 
\end{equation}
where $E_s^i$ has level $n(i)$ and $j(i)<n$ satisfies $n(i) + j(i)\in 
\{j,j+n,j+2n,\ldots\}$. 
Since $\|w\|=\|w_i\|= 1$, $\sum_{i\in A} 2^{-n(i)} |a_i| \le 1$. 
Also 
$$\frac12 \sum_{s=1}^{\ell (i)} |E_s^i w_i|_{j(i)}^n 
\le |w_i|_{j(i)+1}^n \le Lc$$ 
by our hypothesis. 
Hence 
$$|w|_j^n \le 2Lc\ .$$

Of course $\T$ may badly split some $w_i$'s. 
In this case we alter $\T$ as follows. 
Starting with the smallest level we check to see if a given level 
badly splits any $w_i$'s. 
If it does we split the offending sets at $\min\supp w_i$ and 
$\max \supp w_i$. 
Thus a given $E\in \T$ could be split into at most 3 pieces at this stage. 
We intersect successors of split sets with each of the at most three new 
pieces maintaining a tree, but losing admissibility. 
Then proceed to the next level of the new tree and repeat. 
We now have a tree $\T'$ that does not badly split any $w_i$. 
If we replace each set $E$ in this tree by $3E$ we obtain an admissible tree. 
Thus $\T'$ is admissible (3). 
Furthermore, we obtain an expression like \EQ{2.1}, except that the
equality is replaced by the inequality:
$$|w|_j^n \le \sum_{i\in A} 2^{-n(i)} |a_i| \sum_{s=1}^{\ell (i)} 
|E_s^i w_i|_{j(i)}^n\ ,$$ 
where the sets $(E_s^i)$ come from our altered tree just like \EQ{2.1} 
was obtained from $\T$.

Letting $m_i = \min \supp w_i$  we have 
$\|\sum a_i e_{3{m_i}}\|\le K_3 \|\sum a_i e_{m_i}\| \le K_3 \|w\| =
K_3$. Since  $\T'$ is an admissible (3) tree we have 
$$\sum_{i\in A} 2^{-n(i)} |a_i| \le K_3\ .$$ 
Thus $|w|_j^n \le 2K_3 Lc$.
\end{proof}

\pf{Proof of  Theorem~\ref{thm:2.1}}
Fix  $0<\ep <1$ to be specified later. 
By Proposition~\ref{prop:1.1} we may assume 
\begin{equation}\label{eq:2.3} 
\Big|\, |y|_j^n - |y|_{j+n}^n\, \Big|  <\ep\ \mbox{ for }\ 0\le j\le n
\ \mbox{ and }\ y\in S_Y\ .
\end{equation}
Also we may assume $(y_\ell) \prec (e_\ell)$ is a normalized (in $T$) basis 
for $Y$ and that for some $(c_j)_0^{2n}\subseteq (0,1]$, 
\begin{equation}\label{eq:2.4} 
\Big|\, |y_\ell|_j^n - c_j\,\Big| <\ep \ \mbox{ for all }\ \ell\ ,\qquad 
0\le j\le 2n\ .
\end{equation}

Hence we also have, from \EQ{2.3} and \EQ{2.4}, 
\begin{equation}\label{eq:2.5} 
|c_j - c_{j+n}| <3\ep\ \mbox{ if }\ 0\le j\le n\ .
\end{equation}

\begin{lemma}\label{lem:2.3} 
Let $0<i\le n$ and let $z= \sum_{\ell\in A} \alpha_\ell y_\ell$ be an 
$(i,\ep)$ average $(3)$ of $(y_\ell)$, 
$\alpha_\ell >0$ for $\ell\in A$. 
Thus $(y_\ell)_{\ell \in A}$ is $i$ admissible and $\sum_{\ell\in A}
\alpha_\ell=2^i$.
Then 
\begin{eqnarray}
&&c_{j-i} - \ep \le |z|_j^n \le 2c_{j-i+1} + (K_3 +1) 2\ep\ ,\qquad\qquad 
0< i\le j\le n \label{eq:2.6}\\
&&c_{n+j-i} -\ep \le |z|_j^n \le 2c_{n+j-i+1} +(K_3+1)3\ep\ ,\qquad 
0\le j<i\le n \label{eq:2.7}
\end{eqnarray}
\end{lemma}

\begin{proof} 
\EQ{2.6}:  Let $k= j-i$. 
{From} Proposition~\ref{prop:1.1}, \EQ{2.4} and the fact that $S_k[S_i]=S_j$ 
we have 
\begin{eqnarray*} 
|z|_j^n & \ge & \frac1{2^j} \sum_{\ell\in A} \alpha_\ell \sup 
\biggl\{ \sum_{s=1}^r \|E_s y_\ell\|_n :(E_s)_1^r \mbox{is $k$ admissible}\}
\biggr\}\\
&=& \frac1{2^i} \sum_{\ell\in A} \alpha_\ell |y_\ell|_k^n \ge 
\frac1{2^i} \sum_{\ell\in A} \alpha_\ell (c_k -\ep) 
= c_{j-i} -\ep\ .
\end{eqnarray*} 

The second inequality in \EQ{2.6} is more difficult. 
By Proposition~\ref{prop:1.1} there exist $j$ admissible sets $(E_s)_1^r$ 
with 
\begin{equation}\label{eq:2.8} 
|z|_j^n = \frac1{2^j} \sum_{s=1}^r \|E_sz\|_n\ .
\end{equation}
The sets $(E_s)_1^r$ are the terminal sets of an admissible tree $\T$, 
all having level $j$, and we may assume each $E_s\subseteq \bigcup_{\ell\in A}
\supp y_\ell$. 
We adjust the tree $\T$ by splitting some sets if necessary, as  we did in the 
proof of Lemma~\ref{lem:2.2}, to obtain a tree $\T'$ which is admissible~(3) 
and which does not badly split any $y_\ell$, $\ell\in A$. 
It may be that for some $E\in \T'$ we have $E\subseteq \supp y_\ell$ for 
some $\ell$ and level $E<i$. 
We remove all such sets from the tree $\T'$ (replace each $F$ by 
$F\setminus \cup$ such sets and throw out the empty sets thus obtained). 
This gives us a tree $\T''$ which does not badly split any $y_\ell$ 
and for which no set of level $<i$ is contained in $\supp  y_\ell$ for any 
$\ell\in A$. 
$\T''$ is admissible (3). 

Letting $(E'_s)_{s=1}^{r'}$ and $(E''_s)_{s=1}^{r''}$ be the terminal sets 
of $\T'$ and $\T''$, respectively, \EQ{2.8} yields 
\begin{eqnarray} 
|z|_j^n &\le & \frac1{2^j} \sum_{s=1}^{r'} \|E'_s z\|_n \label{eq:2.9}\\
&=& \frac1{2^j} \sum_{s=1}^{r''} \|E''_s z\|_n 
+ \frac1{2^j} \sum_{s\in D} \|E'_sz\|_n \nonumber
\end{eqnarray} 
where $D = \{ 1 \le s\le r': E'_s$ was discarded from $\T'$ in forming
$\T''\}$.
Let 
\begin{eqnarray*} 
B&\equiv & \{\ell \in A: E'_s \subseteq \supp y_\ell\mbox{ for some } 
s\in D\}\\
&=& \{\ell\in A : E\subseteq \supp y_\ell\mbox{ for some } E\in \T' 
\mbox{ with level } E \le i-1\}\ .
\end{eqnarray*} 
Thus $B'\equiv B\setminus \min B$ satisfies $(y_\ell)_{\ell\in B'}$ is 
$i-1$ admissible (3). 
Hence $\sum_{\ell\in B} \alpha_\ell \le \alpha_{\min B} + \sum_{\ell\in B'} 
\alpha_\ell < \ep +\ep =2\ep$. 
Now $(E'_s)_{s\in D}$ is $j$ admissible (3) and so for $\ell\in B$,
$$\frac1{2^j} \sum_{s\in D} \|E'_s y_\ell \|_n \le \|\tilde y_\ell\|
\le K_3$$ 
where $\tilde y = \sum a_i e_{3i}$ if $y= \sum a_i e_i$. 
Thus 
$$\frac1{2^j} \sum_{s\in D} \|E'_s z\|_n \le K_3 \sum_{\ell\in B}
\alpha_\ell < 2K_3 \ep\ .$$
{From} this and \EQ{2.9} we obtain
\begin{equation}\label{eq:2.10} 
|z|_j^n \le \frac1{2^j} \sum_{s=1}^{r''} \|E''_s z\|_n + 2K_3 \ep\ .
\end{equation} 

Recall that $k = j - i$. 
For $\ell\in A$, $\{ E''_s : E''_s \subseteq \supp y_\ell\}$ 
is $k+1$ admissible. 
Indeed $y_\ell$ could first be split into a 1-admissible family only at 
level $i$ or later by $\T''$. 
The tree $\T''$ continues from this point in an admissible fashion up to 
level $j$. 
Thus from \EQ{2.10}, 
\begin{eqnarray*} 
|z|_j^n &\le & \frac1{2^j} \sum_{\ell\in A} \alpha_\ell \sup 
\biggl\{ \sum_{s=1}^p \|F_s y_\ell\|_n : (F_s)_1^p\mbox{ is $k+1$ admissible}
\biggr\} + 2K_3 \ep\\
&\le & \frac1{2^i} \sum_{\ell\in A} \alpha_\ell 2|y_\ell|_{k+1}^n +2K_3\ep\\
&\le & 2c_{k+1} + 2K_3 \ep + 2\ep\ .
\end{eqnarray*}
This completes the proof of \EQ{2.6}. 

For the lower estimate in \EQ{2.7} note that 
\begin{eqnarray}
|z|_j^n &\ge & |z|_{j+n}^n \ge \frac1{2^{j+n}} \sum_{\ell\in A} \alpha_\ell 
\sup \biggl\{\sum_{s=1}^r \|E_sy_\ell\|_n : (E_s)_1^r
\mbox{ is $n+j-i$ admissible}\biggr\}\label{eq:2.7a}\\
&\ge & \frac1{2^i} \sum_{\ell\in A} \alpha_\ell (c_{n+j-i} - \ep) 
= c_{n+j-i} -\ep\ .\nonumber
\end{eqnarray}
Furthermore the argument in proving the upper estimate 
of \EQ{2.6} yields that 
$|z|_{j+n}^n \le 2c_{n+j-i+1} + (K_3 +1) 2\ep$ and since $|z|_j^n 
\le |z|_{j+n}^n +\ep$ we obtain \EQ{2.7}.
\end{proof}

We continue the proof of Theorem~\ref{thm:2.1} by using
Proposition~\ref{prop:1.2} to construct a block basis $(z_i)_{i=1}^n$
of $(y_\ell)$ so that each $z_i$ is an $(i,\ep)$ average (3) of
$(y_\ell)_{\ell=n}^\infty$.
Let $z= \frac1n \sum_1^n z_i$. 
Let $c= \frac1n \sum_1^n c_i$. 

\begin{lemma}\label{lem:2.4} 
For $0\le j\le n$,  
$$\frac12 c - \frac{n+3}{2n} \ep \le |z|_j^n \le 2c + 3\ep (K_3 +1)$$
\end{lemma} 

\begin{proof} 
If $1\le j< n$ then by Lemma~\ref{lem:2.3}, 
\begin{eqnarray*} 
|z|_j^n & \le & \frac1n \sum_{i=1}^n |z_i|_j^n \\
&\le & \frac1n \Bigl[ 2(c_j + c_{j-1} + \cdots + c_1 + c_n + c_{n-1} 
+ \cdots + c_{j+1}) + (K_3 +1) 3n\ep \Big] \\ 
&\le& 2c + 3\ep (K_3+1)\ .
\end{eqnarray*} 
Similarly if $j=0$ or $n$ 
$$|z|_j^n \le \frac1n \Big[ 2(c_n + c_{n-1} + \cdots + c_1) 
+ (K_3 +1) 3n\ep\Big]\ .$$ 
Hence the upper estimate is established. 

To obtain the lower estimate we note that $(z_i)_1^n$ is 1 admissible
hence by Proposition~\ref{prop:1.1}(b) if $1\le j\le n$
\begin{eqnarray*} 
|z|_j^n & = & \frac1n \Big| \sum_{i=1}^n z_i\Big|_j^n 
\ge \frac12\ \frac1n \sum_{i=1}^n |z_i|_{j-1}^n \\
&\ge & \frac1{2n} [ c_0 + \cdots + c_{n-1} -n\ep]\ .
\end{eqnarray*} 
Since $c_n>c_0-3\ep$, 
$$|z|_j^n \ge \frac1{2n} [c_1+\cdots + c_n - (n+3)\ep] = \frac12 c 
- \frac{n+3}{2n} \ep\ .$$
Also $|z|_0^n \ge |z|_n^n$ and so the lemma is proved.
\end{proof}

Note that by Proposition~\ref{prop:1.2} 
$z$ satisfies $\|z\|\le \max_{1\le i\le n} \|z_i\| \le K_1$ 
and $\|z\| \ge \frac1{2n} \sum_1^n \|z_i\| \ge \frac12$. 

Furthermore, for an arbitrary $y \in T$, $\|y\|=1$ implies that
$|y|_j^n \ge 2^{-n}$ for $0\le j\le n$ and thus we could have
chosen $c_j \ge 2^{-n}$ for $0\le j\le n$ and so in particular
$c\ge 2^{-n}$.
Thus we can choose $\ep$ above to obtain (using Lemma~\ref{lem:2.4})
that the element $z$ satisfies
$$\frac13 c \le |z|_j^n \le 3c\ \mbox{ for }\ 0\le j\le n\ .$$ 
These remarks in conjunction with Lemma~\ref{lem:2.2} complete the proof 
of Theorem~\ref{thm:2.1}.
\endpf

Theorem~\ref{thm:2.1} can be restated as saying that there exists an
absolute constant $K$ such that for all $Y\prec T$ and $n\in \nat$
there exists $Z\prec Y$ satisfying
\begin{equation}
  \label{eq:2.11}
  d = \inf_{0 \le j < n}\inf_{z \in {S_Z}} |z|_j^n 
\le \sup_{0 \le j < n}\sup_{z \in {S_Z}} |z|_j^n \le K d\ .
\end{equation}

It is natural to say that $Z\prec T$ is {\it $n$-stable at $d$} if $Z$
satisfies \EQ{2.11}.
Obvious questions then arise.  
How does $d$ depend upon $Z$, how does it depend upon $n$?
Our next result answers these questions. 

\begin{thm}\label{thm:3.1} 
There exists an absolute constant $L$ so that if  $Z\prec Y$ is
$n$-stable at $d$ then $(Kn)^{-1} \le d\le  L n^{-1}$.
\end{thm} 

\begin{proof} 
The lower estimate is relatively easy. 
Let $Z\prec Y$ be $n$ stable at $d$ and let $z\in S_Z$. 
$\|z\|$ is calculated by a tree ultimately yielding 
$\|z\|=1= \sum_{i\in A} 2^{-n(i)} |z(i)|$ as explained previously. 
The sets in the tree are permitted to stop at any level. 
If we gather together those which stop at levels $j, j+n,j+2n,\ldots$ for 
$j=0,1,\ldots,n-1$ we obtain 
$1\le \sum_{j=0}^{n-1} |z|_j^n$.
Hence for some $j<n$, $|z|_j^n\ge \frac1n$, and thus
$d\ge \frac1{Kn}$, by \EQ{2.11}.

Let $y\in Z \cap [(e_m)]_n^\infty$ with $\|y\| =1$ and $y = \sum a_j e_j$.
For $0\le j<n-1$ choose $y_j^* $ in the unit ball $B_{T^*}$ of $T^*$
so that
$$y_j^* (y) = |y|_j^n = \sum_{s\in {A_j}} 2^{-(j+k_j(s)n)} |a_s| \ge d\ .$$ 
We may assume that $A_j \subseteq \supp y$. 
Note that $\sum_{j=0}^{n-1} y_j^* (y)\ge nd$. 
Partition $\bigcup_{i=0}^{n-1} A_i$ into sets $(E_0,\ldots,E_{n-1})$ 
as follows. 
$s\in E_j$ if and only if  for all $i\ne j$ either $s\notin A_i$ or 
$j+k_j(s)n < i +k_i(s)n$. 
Then $(E_j y_j^*)_{j=0}^{n-1}$ is a collection of $n$ disjointly
supported vectors in $B_{T^*}$ all having support contained in
$[(e_m)]_n^\infty$.
Since $T$ and the modified Tsirelson space $T_M$ are naturally
isomorphic (\cite{CS}) there exists an absolute constant $L'$ so
that
$$\Big\| \sum_{j=0}^{n-1} E_j y_j^*\Big\| \le L' \max_{0 \le j < n}
\|E_j y_j^*\|_{T^*} \le L'\ .$$ 
Furthermore 
$$\biggl( \sum_{j=0}^{n-1} E_j y_j^*\biggr) (y) \ge \frac{nd}2\ .$$

Indeed, for $s \in \bigcup_{j =0}^{n-1} E_j $ pick $j_0$ such that $s
\in E_{j_0}$ and denote  by $F_s$ the set of all $0 \le i < n$, $i \ne
j_0$, such that $s \in A_i$.
Then $\left\{ i+k_i(s)n : i \in F_s\right\} $ is a subset of
$\left\{j_0 + k_{j_0}(s)n+1, j_0 + k_{j_0}(s)n+2, \ldots \right\}$.
Thus
\begin{eqnarray*}
  nd &\le & \sum_{j=0}^{n-1} y_j^* (y)  
  =  \sum_{j=0}^{n-1} \sum_{s \in {E_j}}|a_s| 
   \Bigl(2^{-(j + {k_{j}}(s)n)} + \sum_{i \in {F_s}} 
         2^{-(i + {k_i}(s)n)}\Bigr)\\
  &\le & \sum_{j=0}^{n-1} \sum_{s \in {E_j}}|a_s| 
          \Bigl( 2^{-(j + {k_{j}}(s)n)}  +  2^{-(j + {k_{j}}(s)n)} \Bigr)
   =  2\sum_{j=0}^{n-1} E_j y_j^* (y)\ .  
\end{eqnarray*}
Hence ${nd}/2 \le L'$ so $d\le {2L'}/n$.
\end{proof}

As an immediate consequence of Theorems~\ref{thm:2.1} and
\ref{thm:3.1} we get the following.

\begin{cor}\label{cor:2.6a} 
  There exists an absolute constant $C$ so that for every $Y\prec T$
  and $n\in\nat$ there exists $Z\prec Y$ and $d>0$ so that $Z\prec Y$
  is $n$-stable at $d$ and $(Cn)^{-1} \le d\le  C n^{-1}$.
\end{cor}

\section{Further results}\label{sec3}
\setcounter{equation}{0}

We now turn to some stabilization results for more general norms on
$T$.  
Given an arbitrary equivalent norm $|\cdot|$ on $Y\prec T$, we
describe some procedures on $|\cdot|$, natural in the context of
Tsirelson space, which lead to new norms that cannot distort $T$ by
too much.

Recall \cite{OTW} that if $(y_i)$ is a basis for $Y$ and $n\in \nat$ 
then 
\begin{equation}
  \label{eq:3.1a}
\delta_n(y_i) = \inf \biggl\{\delta\ge0: \Bigl\|\sum_1^k x_i\Bigr\| 
\ge\delta \sum_1^k \|x_i\| \mbox{ whenever }(x_i)_1^k \mbox{ is } 
n\mbox{-admissible} \mbox{ w.r.t. }  (y_j) \biggr\}\ .  
\end{equation}
A result of the type we pursue and which we shall need later was proved 
in \cite{OTW}, Theorem~6.2 (in stronger form). 

\begin{prop}\label{prop:3.2} 
There exists $D<\infty$ so that if $(y_i)$ is a normalized block basis of 
$(e_i)$ for $Y\prec T$ and $|\cdot|$ is an equivalent norm on $Y$ with 
$\delta_1 ((y_i),|\cdot|) = \frac12$ then $|\cdot|$ does not $D$ distort $Y$.
\end{prop} 

\begin{remark}\label{rem:3.3} 
It was shown in \cite{OTW} that for a block basis $(y_i)$ of
$(e_i)$ and any equivalent norm $|\cdot|$ on $Y = [(y_i)]$,
$$\delta_n ((y_i),|\cdot|) \le 2^{-n}\ \mbox{ for all }\ n \ .$$
\end{remark} 

If $|\cdot|$ is an equivalent norm on $Y=[(y_j)]\prec T$ we set for $j\ge0$ 
and $x\in Y$, 
$$|x|_j = \frac1{2^j} \sup \biggl\{ \sum_1^\ell |E_i x| : 
(E_i)_1^\ell \mbox{ is $j$ admissible}\biggr\}\ ,$$ 
(If $x= \sum a_i y_i$, $Ex = \sum_{i\in E} a_i y_i$.) 
Thus $|z|_0 = |z|$ and $|\cdot|_j$ is an equivalent norm on $Y$ for all $j$. 
For $n\in \nat$ we let 
$$|z|^{(n)} = \frac1n \sum_{j=0}^{n-1} |z|_j\ .$$

\begin{prop}\label{prop:3.4} 
There exists $D<\infty$ so that if $n\in \nat$ and  $|\cdot|$ is an
equivalent norm on $Y\prec T$ having basis $(y_i)\prec (e_i)$ and
satisfying
$$\Big| \sum_1^k x_i\Big| \ge \frac12 \sum_1^k |x_i|_{n-1}$$ 
for all $1$ admissible $(x_i)_1^k \prec (y_j)$ then 
$|\cdot|^{(n)}$ cannot $D$ distort  $Y$. 
\end{prop}

\begin{proof} 
Let $(x_i)_1^k$ be 1 admissible w.r.t.\ $(y_\ell)$. 
Then for $j\ge 1$, 
$$\Big| \sum_{i=1}^k x_i\Big|_j \ge \frac12 \sum_{i=1}^k |x_i|_{j-1}$$ 
since $S_1 [S_{j-1}] = S_j$. 
Thus using the hypothesis, 
\begin{eqnarray*} 
\Big| \sum_{i=1}^k x_i\Big|^{(n)} 
& = & \frac1n \sum_{j=1}^{n-1} \Big| \sum_{i=1}^k x_i\Big|_j 
+ \frac1n \Big| \sum_{i=1}^k x_i\Big| \\ 
&\ge & \frac12\ \frac1n \sum_{j=0}^{n-2} \sum_{i=1}^k |x_i|_j 
+ \frac12 \ \frac1n \sum_{i=1}^k |x_i|_{n-1} \\ 
&= & \frac12 \sum_{i=1}^k \biggl( \frac1n \sum_{j=0}^{n-1} |x_i|_j\biggr) 
= \frac12 \sum_{i=1}^k |x_i|^{(n)}\ .
\end{eqnarray*} 
Thus $\delta_1 ((y_i),|\cdot|^{(n)}) =\frac12$. 
The proposition follows from Proposition~\ref{prop:3.2}.
\end{proof} 

\begin{remark}\label{rem:3.5}
The hypothesis of Proposition~\ref{prop:3.4} is satisfied if 
$\delta_n(|\cdot|) = 2^{-n}$. 
\end{remark} 

If $|\cdot|$ is an equivalent norm on $Y=[(y_i)]\prec T$,
we define an equivalent norm on $Y$ by 
\begin{eqnarray*} 
|x|_{\Tr} &=& \sup \biggl\{ \sum_{i\in A} 2^{-n(i)} |E_ix| 
: (E_i)_{i\in A} \mbox{ are the terminal sets}\\
&&\mbox{of an admissible tree with level } E_i = n(i)\biggr\}\ . 
\end{eqnarray*} 
Clearly $|\cdot| \le |\cdot|_{\Tr} $  and if $|\cdot| \le \|\cdot\|$
then $|\cdot|_{\Tr} \le \|\cdot\|$. 
Note that $\|\cdot\| = \|\cdot\|_{\Tr}$ if $(y_i)= (e_i)$. 

The constant $K_2$ appearing in several arguments below is the
constant from Proposition~\ref{prop:1.3}.

\begin{prop}\label{prop:3.6} 
There exists $K$ ($= 2 K_2 M$) so that if $|\cdot|$ is any equivalent
norm on $Y= [(y_i)]\prec T$ then $|\cdot|_{\Tr}$ does not $K$ distort
$Y$.
\end{prop} 

\begin{proof} 
By multiplying $|\cdot|$ by a scalar and passing to $Z\prec Y$ we may
assume $\|\cdot\| \ge |\cdot|$ on $Z$ and $Z$ has a basis
$(z_i)_i^\infty$ with $\|z_i\|=1$ and $|z_i| >\frac12$ for all $i$.
Furthermore by \cite{AO} we may assume that for all $j$ if
$(z_i)_{i \in E}$ is $j$-admissible w.r.t.\ $(e_i)$ then $(z_{i+1})_{i
  \in E}$ is $j$-admissible w.r.t.\ $(y_i)$.

Let $z= \sum_1^\ell a_i z_i$ with $\|z\|=1$. 
Then $\|\sum_1^\ell a_i z_{i-1}\| \ge M^{-1}$ for some absolute 
constant $M$  \cite{CS}.
By Proposition~\ref{prop:1.3} there exists an admissible tree w.r.t.\ 
$(e_i)$ having terminal sets of the form $\supp z_{i-1}$ and level $n(i)$ 
for all $i$ in some set $A$ so that 
$$
\sum_{i\in A} 2^{-n(i)} |a_i| \ge K_2^{-1}\|\sum a_i z_{i-1}\| 
\ge  (K_2 M)^{-1}\ .
$$
It follows that 
$$ |z|_{\Tr} \ge  \sum_{i\in A} 2^{-n(i)} |a_i|\, |z_i|
>  (2K_2 M)^{-1}\ , $$
completing the proof
\end{proof}

\begin{remark}\label{rem:3.6a}
It follows from Proposition~\ref{prop:3.6} that if $|\cdot|$ is an
equivalent norm on $Y= [(y_i)]\prec T$ satisfying $ |y|_{\Tr}\le
\gamma |y|$ for all $y \in Y$, then $|\cdot|$ does not $K\gamma$
distort $Y$.
\end{remark} 

\begin{prop}\label{prop:3.7} 
For all $\gamma >0$ there exists $D(\gamma)<\infty$ with the following 
property. 
Let $Y= [(y_i)]\prec T$. 
If $|\cdot|$ is an equivalent norm on $Y$ and $n\in \nat$ is such that 
$\delta_n ((y_i),|\cdot|) = 2^{-n}$ and $|y|_j \ge \gamma |y|$ 
for all $y\in Y$ and $j<n$, then $|\cdot|$ does not $D(\gamma)$ distort $Y$. 
\end{prop}

\begin{proof} 
By Theorem~\ref{thm:2.1} we may choose $(z_i)\prec (y_i)$, $Z= [(z_i)]$ so 
that for some $d>0$, 
$$d\le \|z\|_n \le Kd\mbox{ for all } z\in S_Z\ .$$
Furthermore by passing to a block basis of $Z$ and scaling $|\cdot|$ as 
necessary we may assume that $\|z\|_n \ge |z|$ for all $z\in Z$ and 
$1= \|z_i\| \ge \|z_i\|_n \ge |z_i| \ge \frac12 \|z_i\|_n$ for all $i$. 
Finally again by \cite{AO} we may assume that if $(z_i)_{i\in E}$ is 
$j$-admissible w.r.t.\ $(e_i)$ then $(z_{i+1})_{i \in E}$ is $j$-admissible 
w.r.t.\ $(y_i)$. 

Let $z= \sum a_i z_i$ with $\|z\|=1$. 

As in the proof of Proposition~\ref{prop:3.6} there exists an
admissible tree w.r.t\ $(y_i)$ having terminal sets of the form $\supp
z_i$ and level $n(i)$, $i\in A$, yielding
$$ \sum_{i\in A} 2^{-n(i)} |a_i|\, \|z_i\| \ge (K_2 M)^{-1}\ .$$ 
Choose $0\le j(i) <n$ so that $n(i) + j(i)\in \{0,n,2n,\ldots\}$. 
Since $\delta_n ((y_i),|\cdot|) = 2^{-n}$ we obtain 
\begin{eqnarray*} 
|z| & \ge & \sum_{i\in A} 2^{-n(i)} |a_i|\, |z_i|_{j(i)}
\ge  \gamma \sum_{i\in A} 2^{-n(i)} |a_i|\, |z_i|\\
&\ge & \frac{\gamma}2 \sum_{i\in A} 2^{-n(i)} |a_i|\, \|z_i\|_n 
\ge  \frac{\gamma d}{2K_2 M}\ .
\end{eqnarray*} 
Thus 
$$\|z\|_n \ge |z| \ge \frac{\gamma}{2K_2K M} \|y\|_n\ .$$ 
Hence 
$Kd \ge |z| \ge \frac{\gamma}{2K_2K M}d$.
The theorem is proved with  $D(\gamma)= 2 \gamma^{-1}K_2 K^2M$
\end{proof} 

Our next result combines the proofs of Proposition~\ref{prop:3.7} and 
the main theorem. 

\begin{prop}\label{prop:3.8} 
For $\gamma>0$ there exists $D(\gamma) <\infty$ so that the following holds. 
Let $n\in \nat$ and let $|\cdot|$ be an equivalent norm on $Y=[(y_i)]
\prec T$ with $\delta_n(|\cdot|) =2^{-n}$. 
Suppose that for all $y\in Y$, $|y|_n \ge \gamma |y|$ and 
$|y| \ge \gamma |y|_j$ for $1\le j\le n$. 
Then $|\cdot|$ does not $D(\gamma)$ distort $Y$. 
\end{prop}

\begin{proof} 
As in the proof of Proposition~\ref{prop:3.7} we may assume that 
$\|\cdot\|_n \ge |\cdot|$ on $Z$, $Z$ has a normalized (in $T$) basis 
$(z_i) \prec (y_i)$ with $|z_i| \ge \frac12 \|z_i\|_n$ for all $i$. 
In addition from Theorem~\ref{thm:2.1} we may assume 
$$d\le |z|_j^n \le Kd \mbox{ for } 0\le j\le n \mbox{ and } z\in S_Z\ .$$
Finally we again assume that if $(z_i)_E$ is $j$-admissible w.r.t.\ $(e_i)$ 
then $(z_{i+1})_E$ is $j$-admissible w.r.t.\ $(y_i)$. 

Note that the hypothesis $\delta_n(|\cdot|) = 2^{-n}$ implies $|\cdot| \ge 
|\cdot|_n$ and more generally $|\cdot|_j \ge |\cdot|_{n+j}$. 
$|\cdot|_n \ge \gamma |\cdot|$ implies that (on $Y$) $|\cdot|_j \le 
\gamma^{-1} |\cdot|_{n+j}$. 

Furthermore we may assume that, for a suitably small $\ep>0$, 
$\big|\, |z_\ell|_j - c_j\big| <\ep$ for all $\ell\in \nat$ and 
$0\le j\le n$ for some $(c_j)_0^n \subseteq \real^+$. 

Fix $1\le i\le n$ and let $z= \sum \alpha_\ell z_\ell$ be an 
$(i,\ep)$ average (3) of $(z_\ell)$. 
Note that $|z|_i \ge \frac1{2^i} \sum_{\ell\in A} \alpha_\ell |z_\ell| \ge
\frac{d}2$ hence 
\begin{description}
\item[(i)] $|z| \ge \gamma |z|_i \ge \frac{d\gamma}2$.
\end{description}

The argument of Lemma~\ref{lem:2.3} remains valid for estimates on $|z|_j$. 
The proof of the upper estimate of \EQ{2.7} yields 
\begin{eqnarray*} 
|z|_{j+n} & \le & 2c_{n+j-i+1} +2\ep (K_3+1)\ ,\mbox{ hence}\\ 
|z|_j &\le & \gamma^{-1} \Big( 2c_{n+j-i+1} +2\ep (K_3 +1)\Big)\ .
\end{eqnarray*} 
If we set $w= \frac1n \sum_1^n w_i$ where $(w_i)_1^n \prec (z_\ell)_n^\infty$ 
and each $w_i$ is an $(i,\ep)$ average (3) of $(z_i)$ we obtain as in 
Lemma~\ref{lem:2.4} that (taking $\ep$ suitably small) 
\begin{description}
\item[(ii)] $\frac13 c \le |w|_j \le 3\gamma^{-1} c$ ($0\le j\le n$) where 
$c= \frac1n \sum_1^n c_i$. 
\end{description} 

Also 
$$|w|_1 \ge \frac1{2n} \sum_1^n |w_i| \ge \frac{d\gamma}4$$ 
from (i) and so 
\begin{description}
\item[(iii)] $|w| \ge \gamma |w|_1 \ge \frac{d}4 \gamma^2$.
\end{description} 

{From} ii) and iii) we have $\frac{d}4 \gamma^2 \le 3\gamma^{-1} c$ and so 
$c\ge \frac{d\gamma^3}{12}$. 
Thus from ii) 
\begin{description}
\item[(iv)] $|w|_j \ge \frac{d\gamma^3}{36}$ for $0\le j\le n$.
\end{description}

We are ready to apply the proof of Proposition~\ref{prop:3.7}. 
Let $(w_i)\prec (z_\ell)$ be such that each $w_i$ is constructed 
as  was $w$ above. 
Let $w= \sum a_i w_i$ with $\|w\|=1$. 
Choose an admissible tree having terminal sets $\supp w_i$ for $i\in A$ 
yielding $\sum 2^{-n(i)} |a_i|\, \|w_i\| \ge (K_2 M)^{-1}$. 
It follows that if $0\le j(i)<n$ satisfies $n(i) + j(i)\in 
\{0,n,2n,\ldots\}$  then 
$$Kd  \ge \|w\|_n \ge |w|  \ge  \sum 2^{-n(i)} |a_i|\, |w_i|_{j(i)}
\ge  \frac{d\gamma^3}{36K_2 M} \ .$$
The theorem is proved with $D(\gamma) = \frac{36KK_2 M}{\gamma^3}$. 
\end{proof}

In comparison with \EQ{1.3}, it is of interest to consider the mixed
Tsirelson space (see \cite{AD}) $T((S_k, c_k2^{-k})_{k})$, where
$c_k \uparrow 1$.
We then ask whether it also coincides with $T$, or at least, whether
its norm, $|\cdot|$ say, is an equivalent norm on a subspace of $T$.
The following result gives the positive answer to the latter question.
It also indicates that the answer to the former question probably
depends upon the asymptotic behavior of $(c_k)$.
Finally, it should be compared with Example 5.12 from \cite{OTW}
which implies that if $c_k < \delta <1$ then no subspace of $T((S_k,
c_k2^{-k})_{k})$ is isomorphic to a subspace of $T$.

\begin{prop}\label{prop:3.20} 
There exists a block subspace $X \prec T$ such that $c \|x\| \le |x|
\le \|x\|$ for $x \in X$, where $c>0$ is an absolute constant,
independent of the choice of $c_k \uparrow 1$.
\end{prop} 

\pf{Outline of the  proof}
Clearly, $|x| \le \|x\|$ for all $x \in T$.
Choose $n(i) \uparrow \infty$ such that $\prod_{1}^\infty c_{n(i)} >
\frac12$ and $\sum_1^\infty 2^{-n(i)}< 1/4 K_2$.
Let $m(1) = n(1)$ and inductively choose $m(i) \uparrow \infty$
so that $m(i+1) \ge 2(m(i) + n(i))$ for all $i=1, 2, \ldots$.

Choose $(x_i)\prec T$ to be a block basis of $(e_i)$ such that each
$x_i$ is an $(m(i),1)$ average $(1)$ of $(e_i)$.
In particular, $x_i = \sum_{j \in {F_i}} \alpha_j^i e_j$, where
$\alpha_j^i >0$ for $j \in F_i$, $F_i \in S_{m(i)}$ and $ \sum_{j\in
  {F_i}} \alpha_j^i = 2^{m(i)}$.
It is easy to check that $1 \le \|x_i\|\le 2$.

Let $x  = \sum_1^\ell a_i x_i$ with $\|x\|=1$.
By Proposition~\ref{prop:1.3} there exists an admissible tree 
$\T$ having terminal sets of the form $\supp x_i$ and level $p(i)$,
yielding 
$$\sum_{i \in S} 2^{-p(i)}|a_i|\, \|x_i\| \ge 1/K_2\ ,$$
for some $S \subseteq \{1, \ldots, \ell\}$. 
Set $G = \{i \in S: p(i) \le n(i)\}$. 
Note that if $B = S \setminus G$ then $\sum_B 2^{-p(i)} \le \sum_B
2^{-n(i)} < 1/4 K_2$, and so
$$\sum_{i \in B} 2^{-p(i)}|a_i|\, \|x_i\|< 2/4 K_2 = 1/2K_2.$$
Thus
$$\sum_{i \in G} 2^{-p(i)}|a_i|\, \|x_i\| >  1/2K_2.$$

Prune the tree $\T$ so as to only admit terminal sets of the form
$\supp x_i$ for $i \in G$. 
Extend each of these sets $m(i)$ levels in an admissible fashion,
ending at the singletons which form $\supp x_i$, ultimately obtaining
an admissible tree $\T'$.
Since $\|x_i\| \le 2$, it follows that
$$ \sum_{i \in G} 2^{-p(i)}|a_i| \sum_{j \in {F_i}} 2^{-m(i)} 
\alpha_j^i  >  1/4K_2\ ,$$
which can be rewritten as
$$ \sum_{i \in G} \sum_{j \in {F_i}}2^{-p(i)-m(i)}|a_i|
\alpha_j^i  >  1/4K_2\ .$$

For $i \in G$ all elements in the support of $x_i$ are terminal sets
of $\T'$ having level $j(i) \equiv p(i) + m(i)$.
Note that for $i' \in G$, $i' > i$, the definition of $G$ and the
growth condition on $m(i)$ imply that
$$j(i') - j(i) = p(i') + m(i') - p(i) - m(i) \ge
m(i') - n(i) - m(i) \ge n(i)\ .$$

Let  $G = \{i_1, \ldots, i_s\}$  written in the increasing order.
The admissible tree $\T'$ has terminal sets of level $j(i_1)$ which
together equal the support of $x_{i_1}$, of level $j(i_2)$ which
together equal the support of $x_{i_2}$, and so on.
Also, $j(i_{k+1}) - j(i_k) \ge n(i_k) \ge n(k)$.

By considering  all the sets of $\T'$ of level $j(i_1)$ we obtain
$$|x| \ge c_{j({i_1})} \Biggl( 2^{j({i_1})} |a_{i_1}| \sum_{j\in
    F_{i_1}} \alpha_j^{i_1} + 2^{j({i_1})} \sum_{r=1}^{r(1)} |E_r^{i_1}
  x|\Biggr)\ ,$$
where $(E_r^{i_1})$ are the remaining sets in $\T'$ of level
$j({i_1})$ which are disjoint from the support of $x_{i_1}$.

We iterate this estimate next continuing the sets $(E_r^{i_1})$ to
level $j({i_2})$ and so on. 
Ultimately we obtain
\begin{eqnarray*}
 |x|  &\ge &   c_{j({i_1})} \left( 2^{j({i_1})} |a_{i_1}| \sum_{j\in
    F_{i_1}} \alpha_j^{i_1} + 
    c_{j({i_2}) - j({i_1})} \Biggl( 2^{j({i_2})} |a_{i_2}| \sum_{j\in
       F_{i_2}} \alpha_j^{i_2} \Biggr.\right.\\
     & &\qquad + \left.\Biggl.
c_{j({i_3}) - j({i_2})} \Bigl( 2^{j({i_3})} |a_{i_3}| \sum_{j\in
       F_{i_3}} \alpha_j^{i_3} + \ldots \ \Bigr)\Biggr)\right)\ .
\end{eqnarray*}
Since $j({i_{k+1}}) - j(i_k) \ge n(k)$, this yields that
$$   |x| \ge  \prod_{k=1}^\infty c_{n(k)} \Biggl(\sum_{r=1}^s
    2^{-j({i_r})} |a_{i_r}| \sum_{j\in 
       F_{i_r}} \alpha_j^{i_r} \Biggr)
    \ge \frac12 (1/4K_2) = 1/8K_2 \ ,$$
completing the proof.
\endpf

Until now we considered the Tsirelson space $T\equiv T(S_1,2^{-1})$,
its subspaces and renormings. 
Analogous results also hold for Tsirelson spaces $T_\theta
\equiv T(S_1,\theta)$, where $0 < \theta <1$.  
It should be noted  however, that absolute constants will change
to functions depending on $\theta$ (typically of the form $c
\theta^{-1}$ where $c$ is an absolute constant).
 
In particular, let us recall that the space $T_\theta$ admits a
$\theta^{-1}-\ep$ distorted norm for every $\ep >0$ (the proof is
exactly the same as for $T$).
In this context a distortion property of the renorming
$T(S_n,\theta^n)$ of $T_\theta$  might be also of interest.

\begin{prop}\label{prop:3.9}
Let $n\in \nat$ and $0 < \theta <1$. 
Let  $X = T(S_n,\theta^n)$. Every $Y \prec X$ contains  $Z\prec Y$
such that $Z$ is $\theta^{-1}-\ep$ distortable for every $\ep >0$.
\end{prop}

\pf{Outline of the  proof} 
First note that the modulus $\delta_{m}$ defined in \EQ{3.1a}
staisfies: 
For all $Y \prec X$ and $k \in \nat$, $\delta_{nk}(Y)\le
\theta^{n(k-1)+1}$.
Indeed, let $Y \prec X$ and let $(y_i)$ be a normalized basis in $Y$.
Let $0 < \ep < 1$ and $y = \sum_{i \in A} \alpha_i y_i$ be an $(nk,
\ep)$ average $(1)$, satisfying conditions i) and ii) of
Proposition~\ref{prop:1.2}.
(Observe that these two conditions have a purely combinatorial
character, and their validity does not depend on the underlying Banach
space.)
In particular, $ \sum_{i\in A} \alpha_i = 2^{nk}$.
Then $\|y\| \ge \delta_{nk}(Y) 2^{nk}$.
Iterating the definition of the norm $k-1$ times we obtain
$$\|y\| \le \theta^{n(k-1)}\sum_{j=1}^\ell \|E_j y\| + 
             \sum_{i \in B}\alpha_i\ ,$$
where $i \in B$ if $\supp y_i$ is split by some set in the tree of
sets obtained by iterating the norm definition. 
Thus $B \in S_{n(k-1)} $ and $E_1 < \ldots <E_\ell$ is
$n(k-1)$-admissible, and for $s \le \ell$ and $i \in A$ one has
$E_s \cap \supp y_i = \emptyset$ or $E_j \supseteq \supp y_i$.
Thus
$$\|y\| \le \theta^{n(k-1)} \sum_{i \in A}\alpha_i + \ep
  \le \theta^{n(k-1)} 2^{nk} +\ep\ .$$
Comparing this with the lower estimate for $\|y\|$ yields the
required bound for $\delta_{nk}(Y)$.

The supermultiplicativity property $\delta_{nk}(Y) \ge
(\delta_1(Y))^{nk}$ (\cite{OTW}, Prop. 4.11) and the previous
estimate immediately imply that for all $Y \prec X$, $\delta_1(Y) \le
\theta$.

This in turn implies that for every $Y \prec X$ there exists $Z\prec
Y$ such that for every $\ep>0$ there is $k \in \nat$ satisfying the
following: for all $W \prec Z$ there exist $w_1 < \ldots w_k$ in $W$
such that $\|\sum_{i=1}^k w_i\| =1$ and $\sum_{i=1}^k \|w_i\| \ge
\theta^{-1}- \ep$.
If not then stabilizing suitable quantities for $k=1, 2, \ldots$ by
passing to appropriate subspaces, then using a diagonal argument and
the definition of $S_1$, we would get a subspace $Y'$ with
$\delta_1(Y') >\theta$.

Now, given $\ep >0$, define $|\cdot|$ on $Z$ by
$$ |z| = \sup_{E_1< \ldots < E_k} \sum_{i=1}^k \|E_i z\|\ ,$$
where $E z$ is the projection  with respect to the basis of $Z$.
Clearly, $\|z\| \le |z| \le k \|z\|$ for $z \in Z$. 
Let  $W \prec Z$.
By the previous claim, there exists $w \in W$ with $\|w\|=1$ and $|w|
\ge \theta^{-1} - \ep$.
On the other hand, a standard argument involving long $\ell_1^m$
averages implies that there exists $x \in W$ with $\|x\|=1$ and 
$|x| \le 1+\ep$ (see e.g., \cite{OTW}, Prop. 2.7).
\endpf

\section{Problems}\label{sec4}
\setcounter{equation}{0}

Of course the main problem is 

\begin{problem}\label{prob:4.1}
Is $T$ arbitrarily distortable? 
Is any subspace of $T$ arbitrarily distortable?
\end{problem}

Our work in Section~\ref{sec3} suggests the following problems. 

\begin{problem}\label{prob:4.2} 
Prove that the class of equivalent norms on $T$ for which 
$\delta_n(|\cdot|) = 2^{-n}$ for some $n>1$ do not arbitrarily distort 
$T$ or any $Y\prec T$. 
\end{problem}

\begin{problem}\label{prob:4.3} 
Prove that for $\gamma >0$ there exists $K(\gamma)<\infty$ so that if 
$|\cdot|$ is an equivalent norm on $T$ satisfying $\delta_n(|\cdot|) \ge 
\gamma 2^{-n}$ for all $n$ then $|\cdot|$ does not $K(\gamma)$ distort 
any $Y\prec T$. 
\end{problem} 

\begin{problem}\label{prob:4.4} 
Prove there exists $K<\infty$ so that if $|\cdot|$ is an equivalent norm 
on $T$ and $Y\prec T$ then for some $n$, $|\cdot|^{(n)}$ does not $K$ 
distort $Y$.
\end{problem}

\medskip\noindent
Edward Odell: 
Department of Mathematics,
The University of Texas at Austin,\\
Austin, TX 78712, USA \\
e-mail: {\small\tt% 
  odell@math.utexas.edu}

\medskip\noindent
Nicole Tomczak-Jaegermann: 
Department of Mathematical Sciences,
University of Alberta,\\
Edmonton, Alberta, T6G 2G1, Canada\\
e-mail: {\small\tt% 
  ntomczak@math.ualberta.ca}

\end{document}